\documentclass[a4paper,11pt]{amsart}
\usepackage{amssymb}
\usepackage[french]{babel}
\usepackage[T1]{fontenc}

\theoremstyle{plain}
\newtheorem{theorem}[subsection]{Th\'{e}or\`{e}me}

\newtheorem{lemma}[subsection]{Lemme}
\newtheorem{proposition}[subsection]{Proposition}

\theoremstyle{definition}

\theoremstyle{remark}

\numberwithin{equation}{section}

\newcommand{\ki}{{\mathcal I}}

\newcommand{\ko}{{\mathcal O}}

\newcommand{\ku}{{\mathcal U}}

\newcommand{\kx}{{\mathcal X}}

\newcommand{\IC}{{\mathbb C}}

\newcommand{\IP}{{\mathbb P}}

\newcommand{\IZ}{{\mathbb Z}}
\newcommand{\gotha}{{\mathfrak a}}
\newcommand{\gothb}{{\mathfrak b}}

\newcommand{\gothd}{{\mathfrak d}}

\newcommand{\gothg}{{\mathfrak g}}

\newcommand{\gothm}{{\mathfrak m}}
\newcommand{\gotho}{{\mathfrak o}}

\newcommand{\gothq}{{\mathfrak q}}

\DeclareMathOperator{\id}{id}
\DeclareMathOperator{\ad}{ad}

\newcommand{\lra}{\longrightarrow}
\newcommand{\xra}{\xrightarrow}
\newcommand{\isom}{\cong}
\newcommand{\tensor}{\otimes}

\newcommand{\cf}{\emph{cf.}\ }

\DeclareMathOperator{\Sym}{Sym}
\DeclareMathOperator{\Hom}{Hom}
\DeclareMathOperator{\End}{End}
\DeclareMathOperator{\Ext}{Ext}
\DeclareMathOperator{\Aut}{Aut}
\DeclareMathOperator{\Bild}{Im}
\DeclareMathOperator{\Kern}{Ker}
\DeclareMathOperator{\Kokern}{Coker}

\DeclareMathOperator{\Liesp}{\mathfrak{sp}}
\DeclareMathOperator{\LieSp}{Sp}
\DeclareMathOperator{\Liesl}{\mathfrak{sl}}
\DeclareMathOperator{\LieSl}{SL}
\DeclareMathOperator{\Liegl}{\mathfrak{gl}}
\DeclareMathOperator{\LieGl}{GL}

\DeclareMathOperator{\Grass}{Grass}
\DeclareMathOperator{\Proj}{Proj}
\DeclareMathOperator{\Spec}{Spec}

\DeclareMathOperator{\Bl}{Bl}
\DeclareMathOperator{\Def}{Def}
\DeclareMathOperator{\sing}{sing}

\DeclareMathOperator{\pair}{pair}

\DeclareMathOperator{\ch}{ch}
\DeclareMathOperator{\trace}{tr}

\DeclareMathOperator{\ini}{in}
\DeclareMathOperator{\gr}{gr}

\newcommand{\bull}{\bullet}
\newcommand{\tZ}{\widetilde{Z}}
\newcommand{\oB}{\overline{B}}
\newcommand{\ob}{\overline{b}}
\newcommand{\smb}{{\scriptstyle\bullet}}
\newcommand{\tdph}{termes de degree plus haut}
\newcommand{\kit}{J}

\begin{document}

\title{La singularit\'{e} de O'Grady}
%\date{22 Septembre 2005}

\author{Manfred Lehn \& Christoph Sorger}
\address{
Manfred Lehn\\
Institut f\"{u}r Ma\-the\-ma\-tik\\
Johannes Gu\-ten\-berg-Uni\-ver\-si\-t\"{a}t Mainz\\
D-55099 Mainz, Allemagne } \email{lehn@mathematik.uni-mainz.de}

%\author{Christoph Sorger}
\address{
Christoph Sorger\\
Institut Universitaire de France \& Laboratoire de Math\'{e}matiques Jean Leray (UMR 6629 du CNRS)\\
Universit\'{e} de Nantes\\
2, Rue de la Houssini\`{e}re\\
BP 92208\\
F-44322 Nantes Cedex 03, France}
\email{christoph.sorger@univ-nantes.fr}


\begin{abstract}
Soit $M_{2v}$ l'espace de modules des faisceaux semi-stables de
vecteur de Mukai $2v$ sur une surface K3 ou ab\'{e}lienne o\`{u} $v$ est
primitif tel que $\langle v,v \rangle=2$. Nous montrons que
l'\'{e}clatement de $M_{2v}$ le long de son lieu singulier r\'{e}duit
fournit une r\'{e}solution symplectique des singularit\'{e}s. Ceci donne une
description directe des r\'{e}solutions de O'Grady de $M_{K3}(2,0,4)$ et
$M_{Ab}(2,0,2)$.

\noindent{\textsc{Abstract.}} Let  $M_{2v}$ be the moduli space of
semistable sheaves with \nobreak{Mukai} vector $2v$ on an abelian or
K3 surface where $v$ is primitive such that $\langle v,v \rangle=2$.
We show that the blow-up of the reduced singular locus of $M_{2v}$
provides a symplectic resolution of singularities. This provides a
direct description of O'Grady's resolutions of $M_{K3}(2,0,4)$ and
$M_{Ab}(2,0,2)$.
\end{abstract}

\maketitle

%%%%%%%%%%%%%%%%%%%%%%%%%%%%%%%%%%%%%%%%%%%%%%%%%%%%%%%%%%%%%%%%%%%%%%%%%%%%%
%%%%%%%%%%%%%%%%%%%%%%%%%%%%%%%%%%%%%%%%%%%%%%%%%%%%%%%%%%%%%%%%%%%%%%%%%%%%%
\section{L'espace de modules de O'Grady}

Il est notoirement difficile de produire des vari\'{e}t\'{e}s
holomorphiquement symplectiques irr\'{e}ductibles compactes. \`{A}
d\'{e}formation pr\`{e}s, il y a les deux s\'{e}ries infinies trouv\'{e}es par
A.~Beauville \cite{B}, \`{a} savoir les sch\'{e}mas de Hilbert des points
sur une surface $K3$ et les vari\'{e}t\'{e}s de Kummer g\'{e}n\'{e}ralis\'{e}es
associ\'{e}es \`{a} une surface ab\'{e}lienne, et les deux exemples isol\'{e}s
construits par O'Grady (\cite{OG1}, \cite{OG2}).

Selon Mukai, les espaces de modules de faisceaux stables sur une
surface $K3$ ou ab\'{e}lienne sont lisses et munis d'une forme
symplectique holomorphe. Si les param\`{e}tres num\'{e}riques, \`{a} savoir le
rang et les classes de Chern, sont choisis de fa\c{c}on \`{a} ce qu'il
n'existe pas de faisceau strictement semi-stable, les espaces de
Mukai sont compacts mais se d\'{e}forment dans les exemples de
Beauville.

La strat\'{e}gie de O'Grady pour la construction de ses exemples est de
consid\'{e}rer certains espaces de modules singuliers puis de les
r\'{e}soudre symplectiquement si possible. La r\'{e}solution de O'Grady est
relativement compliqu\'{e}e: elle consiste en deux \'{e}clatements suivis
d'une contraction sur un sch\'{e}ma de param\`{e}tres muni d'une action d'un
groupe de reparam\'{e}trisation o\`{u} il faut tenir compte \`{a} chaque \'{e}tape
du comportement de la semi-stabilit\'{e} des points sous-jacents.

Nous montrons ici que l'\'{e}clatement du lieu singulier r\'{e}duit de
l'espace de modules fournit directement une r\'{e}solution symplectique
des singularit\'{e}s. La m\'{e}thode consiste en une \'{e}tude locale d\'{e}taill\'{e}e
des singularit\'{e}s de l'espace de module sans r\'{e}f\'{e}rence aux techniques
de constructions globales.

Soit $X$ une surface projective lisse de type $K3$ ou ab\'{e}lien. Soit
de plus $v\in H^{\pair}(X,\IZ)$ et soit $M_v$ l'espace de modules
des faisceaux semi-stables sur $X$ par rapport \`{a} un diviseur ample
$H$ et ayant vecteur de Mukai $v$ (voir section \ref{sec:conenormal}).

\begin{theorem}\label{th:main}--- Soit $v\in H^{\pair}(X,\IZ)$ un
\'{e}l\'{e}ment primitif avec $\langle v,v \rangle=2$ et $H$ un diviseur
ample $2v$-g\'{e}n\'{e}rique. Alors l'\'{e}clatement de $M_{2v}$ le long de
son lieu singulier r\'{e}duit est une r\'{e}solution symplectique des
singularit\'{e}s.
\end{theorem}

Les hypoth\`{e}ses du th\'{e}or\`{e}me couvrent \`{a} la fois les exemples de
O'Grady \cite{OG1,OG2} et aussi celles, similaires, de Rapagnetta
\cite{Rapagnetta}. Un ph\'{e}nom\`{e}ne similaire a \'{e}t\'{e} observ\'{e} par Haiman
pour la puissance sym\'{e}trique de $\IC^2$ et la r\'{e}solution
symplectique donn\'{e} par le sch\'{e}ma de Hilbert \cite{H1},\cite{H2}. Par
ailleurs, dans \cite{KLS}, nous avons montr\'{e} qu'\`{a} part dans les cas
ci-dessus du th\'{e}or\`{e}me 1.1, les espaces de modules des faisceaux
semi-stables sans torsion sur une surface $K3$ ou ab\'{e}lienne
n'admettent pas de r\'{e}solution symplectique en rang $r\geq 2$ (voir
\cite{KLS}, Th\'{e}or\`{e}me B pour l'\'{e}nonc\'{e} pr\'{e}cis).

Le plan de la d\'{e}monstration est le suivant. Nous allons d'abord
\'{e}tudier la g\'{e}om\'{e}trie d'une certaine orbite nilpotente
$Z\subset\Liesp_4$ qui nous servira comme mod\`{e}le locale pour la plus
mauvaise singularit\'{e} de $M_{2v}$, puis on \'{e}tablira une version
analogue du th\'{e}or\`{e}me \ref{th:main} pour $Z$ (Th\'{e}or\`{e}me
\ref{thm:laresolution}). Ensuite, nous montrerons un r\'{e}sultat
essentiel pour notre preuve qui est aussi d'int\'{e}r\^{e}t ind\'{e}pendant, \`{a}
savoir que $Z$ poss\`{e}de une certaine rigidit\'{e} relative aux
d\'{e}formations: une d\'{e}formation de $Z$ qui ne change pas les
singularit\'{e}s de $Z$ en dehors de l'origine ne peut d\'{e}former la
singularit\'{e} \`{a} l'origine non plus. Plus pr\'{e}cis\'{e}ment, on montrera que
le module $T^1_Z$ est pur (Th\'{e}or\`{e}me \ref{thm:purete}).

Apr\`{e}s ces pr\'{e}parations nous \'{e}tudierons la structure locale de
l'espace de modules. Les singularit\'{e}s de l'espace de modules
proviennent \`{a} la fois du passage au quotient par l'action des
groupes d'automorphismes et du comportement de l'application de
Kuranishi dont nous rappelons la construction et ses propri\'{e}t\'{e}s dans
un appendice (Appendice A). Nous nous placerons ensuite en la plus
mauvaise singularit\'{e} et montrerons que le c\^{o}ne normal de cette
singularit\'{e} est isomorphe au mod\`{e}le affine ci-dessus (Th\'{e}or\`{e}me
\ref{thm:isomorphisme}). Finalement, nous montrerons, par un
argument de d\'{e}formation explicite, que la singularit\'{e} est
formellement isomorphe \`{a} son c\^{o}ne normal.

On voit ainsi que l'\'{e}clatement du lieu singulier r\'{e}duit donne une
r\'{e}solution semi-petite des singularit\'{e}s de $M_{2v}$. Il est connu
que la forme symplectique s'\'{e}tend sur l'image r\'{e}ciproque de la
partie g\'{e}n\'{e}rique du lieu singulier qui consiste en des
singularit\'{e}s de type $A_{1}$ (\cf\cite{OG1}). Cette forme s'\'{e}tend
partout pour des raisons de codimension.

\emph{Remerciements:} Nous remercions Duco van Straten pour le aide
avec les d\'{e}formations ainsi que Theo de Jong et Dmitry Kaledin pour
nos discussions tr\`{e}s constructives. Nous avons fait des exp\'{e}riences
avec \textsc{Singular} \cite{Greuel} et nous remercions Gert-Martin
Greuel pour son aide pendant nos premiers pas avec ce logiciel.
Finalement, nous remercions le referee pour nous avoir signal\'{e} que
notre argument dans la premi\`{e}re version de cet article \'{e}tait
incomplet.

%%%%%%%%%%%%%%%%%%%%%%%%%%%%%%%%%%%%%%%%%%%%%%%%%%%%%%%%%%%%%%%%%%%%%%%%%%%%%%%
%%%%%%%%%%%%%%%%%%%%%%%%%%%%%%%%%%%%%%%%%%%%%%%%%%%%%%%%%%%%%%%%%%%%%%%%%%%%%%%

\section{Le mod\`{e}le alg\'{e}brique et sa r\'{e}solution}
\label{section:resolution}

Soit $V$ un espace vectoriel complexe de dimension 4 muni d'une
forme symplectique $\omega$. On d\'{e}signe par $\Liesp(V)$ l'alg\`{e}bre de
Lie symplectique associ\'{e}. Soit $Z\subset\Liesp(V)$ la sous-vari\'{e}t\'{e}
des $B\in\Liesp(V)$ tels que $B^2=0$. Les 16 coefficients de la
matrice $B^2$ ne sont pas ind\'{e}pendants: l'id\'{e}al $I_0\subset
\IC[\Liesp(V)]$ de $Z$ est engendr\'{e} par 6 quadriques. Il est connue
que $Z$ est l'adh\'{e}rence de l'orbite nilpotente de type $\gotho(2,2)$
(voir \cite{Coll}, \cite{Fu} pour plus de d\'{e}tails). Le lieu
singulier $Z_{\sing}$ de $Z$ consiste en les $B$ ayant rang $\leq
1$; son id\'{e}al $L_0$ est engendr\'{e} par les coefficients de
$\Lambda^2B$, c'est-\`{a}-dire les $2\times 2$-mineurs de $B$.

On sait qu'on obtient une r\'{e}solution de $Z$ de la mani\`{e}re
suivante. Soit $G$ la grassmannienne des sous-espaces vectoriels
$U\subset V$ isotropes maximaux. La vari\'{e}t\'{e} $G$ est isomorphe
\`{a} la sous-vari\'{e}t\'{e} de $\IP(\Lambda^2V^*)=\IP^5$ d\'{e}finie
par l'\'{e}quation quadratique $\varphi\wedge\varphi=0$ et
l'\'{e}quation lin\'{e}aire $\omega(\varphi)=0$ avec
$\varphi\in\Lambda^2V$. Soit $\ku\subset V\otimes\ko_G$ le
sous-fibr\'{e} tautologique. Comme $\omega$ est
non-d\'{e}g\'{e}n\'{e}r\'{e}e et les fibres de $\ku$ sont isotropes par
construction, $\omega$ induit un accouplement parfait
$\ku\times (V\otimes_\IC\ko_G)/\ku\to\ko_G$. On obtient ainsi la suite
exacte tautologique
\begin{equation}
0\lra\ku\lra V\otimes\ko_G\lra\ku^*\lra0.
\end{equation}
Le fibr\'{e} vectoriel $\Hom(\ku^*,\ku)$ est la restriction \`{a} $G$ du
fibr\'{e} cotangent de $\Grass(2,V)$ ; le fibr\'{e} cotangent $T^*G$ de
$G$ s'identifie avec le sous-fibr\'{e} des formes sym\'{e}triques sur
$\ku^*$.

Soit $\tZ\subset Z\times G$ la sous-vari\'{e}t\'{e} des paires
$(B,U)$ telles que $B(U)=0$. D\'{e}signons par $\sigma:\tZ\to Z$ et
$\pi:\tZ\to G$ les deux projections canoniques. Par construction,
pour un \'{e}l\'{e}ment $(B,U)\in\tZ$, l'endomorphisme $B:V\to V$ se
factorise comme suit
\begin{equation}
V\lra U^*\xra{\;\oB\;}U\lra V
\end{equation}
o\`{u} $\oB$ est une forme quadratique sur $U^*$. On voit ainsi que
$\pi:\tZ\to G$ est isomorphe \`{a} la projection canonique $T^*G\to
G$. Par cons\'{e}quent, $\tZ$ est lisse.

Le morphisme $\sigma:\tZ\to Z$ est une r\'{e}solution semi-petite~: la
fibre de $\sigma$ au-dessus d'un point $B\in Z$ de rang $1$ est
isomorphe \`{a} $\IP((\Kern B/\Bild B)^*)=\IP^1$ et au-dessus de $B=0$
isomorphe \`{a} $G$.

\begin{theorem}\label{thm:laresolution}---
La r\'{e}solution $\sigma:\tZ\to Z$ est isomorphe \`{a} l'\'{e}clatement de
$Z$ le long de $Z_{\text{sing}}\subset Z$.
\end{theorem}

\begin{proof} Soit $b:V\otimes\ko_Z\to V\otimes\ko_Z$ l'homomorphisme
canonique qui en $B\in Z$ multiplie par $B$. Par construction $\sigma^*b$
s'annule sur $\pi^*\ku$. On obtient la factorisation
\begin{equation}
\sigma^*b:V\otimes\ko_{\tZ}\lra\pi^*\ku^*\xra{\;\ob\;}\pi^*\ku\lra
V\otimes\ko_{\tZ}
\end{equation}
En passant \`{a} la puissance ext\'{e}rieure
\begin{equation}
\sigma^*\Lambda^2b:\Lambda^2V\otimes\ko_{\tZ}\lra\pi^*\det\ku^*
\xra{\;\Lambda^2\ob\;}\pi^*\det\ku\lra \Lambda^2V\otimes\ko_{\tZ}
\end{equation}
on voit que l'annulation de la section $\Lambda^2\ob$ du fibr\'{e}
inversible $\pi^*(\det\ku)^2$ \'{e}quivaut \`{a} l'annulation de
$\sigma^*\Lambda^2b$. Par cons\'{e}quent,
% si $\widetilde{J}\subset\ko_{\tZ}$ d\'{e}signe le faisceau d'id\'{e}aux de
% $\sigma^{-1}(Z_{\text{sing}})$
on a
\begin{equation}\label{eq:sigma-ample}
J:=\sigma^{-1}L_0\cdot\ko_{\tZ}\isom\pi^*(\det\ku^*)^2\isom
\pi^*\ko_{\IP(\Lambda^2V^*)}(2)
\end{equation}
On observe par ailleurs que $\sigma_*(\ko_{\tZ})=\ko_Z$ puisque
les singularit\'{e}s de $Z$ sont rationnelles. Consid\'{e}rons maintenant
les inclusions
\begin{equation}
L_0\ko_Z\subset\sigma_*J\subset\ko_Z
\end{equation}
Les radicaux de $L_0\ko_Z$ et $\sigma_*J$ co\"{\i}ncident.  Comme
$L_0$ est un id\'{e}al premier, $L_0\ko_Z=\sigma_*{J}$. Selon le lemme
suivant, on a \'{e}galement $L_0^n\ko_Z=\sigma_*({J}^n)$. Ainsi
\begin{equation}%
\tZ\isom\Proj(\bigoplus\sigma_*({J}^n))=\Proj(\bigoplus
L_0^n\ko_Z)=\Bl(Z).
\end{equation}
\end{proof}

\begin{lemma}--- L'homomorphisme $(\sigma_*{J})^{\otimes n}\to
\sigma_*({J}^n)$ est surjectif.
\end{lemma}

\begin{proof}
On consid\`{e}re le diagramme commutatif
$$\begin{array}{ccc}
T^*G & \xra{\pi} & G\\
\scriptstyle{\sigma}\Big\downarrow && \Big\downarrow\\
Z& \lra & *\\
\end{array}
$$
Observons maintenant que l'on a
\begin{equation}
H^0(Z,\sigma_*({J}^{n}))=H^0(G,\pi_*\pi^*\ko(2n))=H^0(G,S^{*}(T_G)\otimes\ko(2n))
.\end{equation}
Comme $Z$ est affine, il suffit de montrer que
$$H^0(G,S^*(T_G)\otimes\ko(2))^{\tensor n}\to
H^0(G,S^{*}(T_G)\otimes\ko(2n))$$ est surjectif. Pour cela, il
suffit de voir que
$$
H^0(G,\ko(2))\tensor H^0(G,S^{k}(T_G)\otimes\ko(2n))\to
 H^0(G,S^{k}(T_G)\otimes\ko(2n+2))
$$
est surjectif pour tout $k\geq 0$ et tout $n\geq 1$ ou encore que
$S^{k}(T_G)$ est $2$-r\'{e}gulier sur $G$ au sens de Mumford-Castelnuovo
\cite{Kleiman}. Que $H^1(G,S^k(T_G)(1))=0$ est cons\'{e}quence imm\'{e}diate
du th\'{e}or\`{e}me d'annulation de Griffiths \cite{Griffiths}. L'annulation
de $H^2(G,S^k(T_G))$ et de $H^3(G,S^k(T_G)(-1))$ se d\'{e}duit des
suites exactes
$$
0 \lra S^{k-1}(\ko_G(1)^{\oplus 5}) \lra S^{k}(\ko_G(1)^{\oplus 5})
\lra  S^k(T_{\IP^4}|_G)\lra 0 $$ et
$$
 0 \lra  S^{k}(T_G)\lra
S^k(T_{\IP^4}|_G) \lra S^{k-1}(T_{\IP^4}|_G)\otimes\ko_G(2)\lra 0.
$$
\end{proof}

%%%%%%%%%%%%%%%%%%%%%%%%%%%%%%%%%%%%%%%%%%%%%%%%%%%%%%%%%%%%%%%%%%%%%%%%%%%%%%
%%%%%%%%%%%%%%%%%%%%%%%%%%%%%%%%%%%%%%%%%%%%%%%%%%%%%%%%%%%%%%%%%%%%%%%%%%%%%%

\section{Rappel sur des d\'{e}formations et Rigidit\'{e}}
\label{section:deformation}

%Commen\c{c}ons par un rappel sur les d\'{e}formations:
Soit $(X,x_0)$
un germe d'un espace analytique. Une d\'{e}formation de $X$ est une
application plate $f:(\kx,x_0)\to (S,s_0)$ des germes d'espaces analytiques
avec un isomorphisme explicite $\kx_0\isom X$ de la fibre sp\'{e}ciale
$\kx_0=f^{-1}(s_0)$. Soit $\Def(X)$ le foncteur qui associe \`{a} une
base $S$ l'ensemble des d\'{e}formations de $X$ sur $S$ modulo
isomorphisme. L'espace tangent formel $\Def(X)(\IC[\varepsilon])$
est le germe d'un faisceau de $\ko_X$-modules coh\'{e}rent $T^1_X$
d\'{e}fini comme suit: Soit $X\subset (\IC^n,0)$ donn\'{e} par
l'id\'{e}al $I\subset \IC\{x_1,\ldots,x_n\}$ et soit
$\Theta:\ko_X\langle \partial_{x_1},\ldots,\partial_{x_n}\rangle\to
\Hom_{\ko_X}(I/I^2,\ko_X)$ l'homomorphisme naturel avec
$\partial_{x_i}\mapsto (f\mapsto \partial_{x_i}(f))$. Alors
$T^1_X=\Kokern(\Theta)$. Par le crit\`{e}re de Jacobi $T_X^1$ a pour support
le lieu singulier de $X$.

Pour une hypersurface $X=\{f=0\}\subset\IC^n$ le calcul de $T^1_X$
est simple:
$T_X^1=\IC\{x_1,\ldots,x_n\}/(f,\partial_{x_1}f,\ldots,\partial_{x_n}f)$.
Par exemple, pour la singularit\'{e} de type $A_1$,
$\{x_1^2+x_2^2+x_3^2=0\}$ on trouve $T_X^1=\IC$.

Le th\'{e}or\`{e}me suivant est essentiel pour la d\'{e}monstration de notre
r\'{e}sultat principal. Il montre une certaine rigidit\'{e} inattendue de la
singularit\'{e} de $Z$:

% Une d\'{e}formation de $Z$ qui ne change pas les singularit\'{e}s
% de $Z$ en dehors de l'origine ne peut d\'{e}former la singularit\'{e} \`{a} l'origine
% non plus.

\begin{theorem}\label{thm:purete}---
Soit $Z\subset\Liesp(V)$ la vari\'{e}t\'{e} d\'{e}finie ci-dessus.
%dans la section \ref{section:resolution}.
Alors le $\ko_Z$--module $T^1_Z$ est pur de dimension $4$.
\end{theorem}

\begin{proof} Rappelons que $Z$ est singulier le long du lieu
de matrices de rang $\leq 1$ dans $\Liesp_4$. En un point $p\in
Z_{sing}\setminus\{0\}$, le germe $(Z,p)$ est essentiellement une
singularit\'{e} de type $A_1$. Plus pr\'{e}cis\'{e}ment, il est analytiquement
isomorphe au germe $(\IC^4,0)\times(\{x_1^2+x_2^2+x_3^2=0\},0)$. On
en d\'{e}duit que la restriction de $T^1_Z$ \`{a} $Z\setminus\{0\}$ est
localement libre de rang $1$ le long de $Z_{sing}\setminus\{0\}$. Un
sous-module de $T^1_Z$ de dimension $<4$ est donc n\'{e}cessairement
concentr\'{e} \`{a} l'origine.

Soit $\ki_0=I_0\ko_{\Liesp(V)}\subset\ko_{\Liesp(V)}$ le faisceau
d'id\'{e}aux de $Z\subset\Liesp(V)$. Comme ci-dessus on d\'{e}signe par
$\Theta:T_{\Liesp(V)}|_Z\rightarrow\Hom(\ki_0,\ko_Z)$ le morphisme
canonique qui associe \`{a} un champ de vecteurs la d\'{e}rivation
correspondante.
% En chaque point $p\in Z$, le germe du conoyau
% $T^1_Z=\Kokern(\Phi)$ en $p$ est isomorphe \`{a} l'espace tangent des
% d\'{e}formations du germe $(Z,p)$.

Une pr\'{e}sentation de l'id\'{e}al $\ki_0$ se construit comme suit.
La structure symplectique $\omega$ sur $V$ induit une
anti-involution $f\mapsto f^t$ sur $\End(V)$ donn\'{e}e par
$\omega(f^t(v),w)=\omega(v,f(w))$ pour $v,w\in V$. Soit
$\End(V)=E_+\oplus E_-$ la d\'{e}composition par rapport aux valeurs
propres de l'involution telle qu'en particulier $\Liesp(V)=E_-$. On a
une surjection
\begin{equation}\label{eq:surjection-sur-I}
  g:E_+^*\otimes\ko_{\Liesp(V)}\lra \ki_0,\quad
  (\lambda,B)\mapsto\lambda(B^2).
\end{equation}
En effet, $Z\subset\Liesp(V)$ est d\'{e}fini par l'\'{e}quation $B^2=0$,
et $B^2\in E_+$ pour $B\in\Liesp(V)$. Il est facile \`{a} v\'{e}rifier que
la suite suivante est un complexe
\begin{equation}\label{eq:presentation-de-I}
  (\Lambda^2 E_+^*\oplus E_+^*)\otimes\ko_{\Liesp(V)}\xra{\ r\ }
   E_+^*\otimes\ko_{\Liesp(V)}\xra{\ g\ } \ki_0,
\end{equation}
o\`{u} l'application $r$ est donn\'{e}e sur la partie $\Lambda^2 E_+^*$
par les relations tautologiques
$r(\lambda\wedge\mu,B)=\lambda\otimes\mu(B^2)-\mu\otimes\lambda(B^2)$
et sur la partie $E_+^*$ par $r(\lambda,B)=\lambda\circ\ad B$. On
peut montrer que \eqref{eq:presentation-de-I} est une pr\'{e}sentation
de $\ki_0$, mais nous n'en avons pas besoin. En appliquant le foncteur
$\Hom(-,\ko_Z)$ on obtient le complexe
\begin{equation}\label{eq:erster-komplex}
  \Hom(\ki_0,\ko_Z)\xra{\ g^\vee\ }
  E_+\otimes\ko_Z\xra{\ r^\vee\ }
  E_+\otimes\ko_{Z}.
\end{equation}
Si $d^0$ d\'{e}signe la composition
\begin{equation}
T_{\Liesp(V)}\otimes\ko_Z=E_-\otimes\ko_Z\xra{\ \Phi\ }
\Hom(\ki_0,\ko_Z)\lra E_+\otimes\ko_Z
\end{equation}
on voit que l'on a  le lemme suivant.
\begin{lemma}\label{lemma:H1}
Le module $T^1_Z$ s'identifie \`{a} un sous-module du $H^1$ du complexe
\begin{equation}
C^\bull:E_-\otimes\ko_Z\xra{\ d^0\ } E_+\otimes\ko_Z \xra{\ d^1\ }
E_+\otimes\ko_Z
\end{equation}
avec $d^0:(b,B)\mapsto bB+Bb$ et $d^1:(b,B)\mapsto bB-Bb$.
\end{lemma}

\noindent Pour d\'{e}montrer le th\'{e}or\`{e}me \ref{thm:purete},
on va se placer sur la r\'{e}solution $\sigma:\widetilde{Z}\rightarrow Z$ de
la section pr\'{e}c\'{e}dente. %\ref{section:resolution}.

\begin{lemma}\label{lemma:H1oben}--- On a
$H^1(C^\bull)=\sigma_*H^1(\sigma^*C^\bull)$.
\end{lemma}

Gr\^{a}ce \`{a} l'\'{e}quivariance de $\sigma^*C^\bull$ sous l'action de
$\LieSp(V)$ les noyaux, conoyaux et groupes de
cohomologie sont plats sur $G$ (par rapport \`{a}
$\pi:\widetilde{Z}=T^*G\to G$). On peut donc se restreindre \`{a}
\'{e}tudier la fibre au-dessus d'un $[U]\in G$.

Choisissons une d\'{e}composition $V=U\oplus U^*$ avec
$\omega=\big(\begin{smallmatrix}\phantom{-}0&|&1\,\\\hline
-1&|&0\,\end{smallmatrix}\big)$. Un \'{e}l\'{e}ment de $E_-$ est de la
forme $\big(\begin{smallmatrix}\alpha &|&\beta \,\\\hline \gamma
&|&-\alpha^*\,\end{smallmatrix}\big)$ avec des $2\times
2$--matrices $\alpha,\beta,\gamma$ o\`{u} $\beta$ et $\gamma$ sont
sym\'{e}triques. De m\^{e}me, un \'{e}l\'{e}ment de $E_+$ est de la forme
$\big(\begin{smallmatrix}\alpha &|&\beta \,\\\hline \gamma
&|&\alpha^*\,\end{smallmatrix}\big)$ avec des $2\times
2$--matrices $\alpha,\beta,\gamma$ o\`{u} $\beta$ et $\gamma$ sont
anti-sym\'{e}triques. La fibre $F$ au-dessus de $U$ est un $\IC^3$
avec coordonn\'{e}es $x,y,z$ o\`{u} $B=\sigma(x,y,z)$ est donn\'{e}e par la
matrice $\big(\begin{smallmatrix}0&|&\overline{B}\\\hline
0&|&0\end{smallmatrix}\big)$ avec $\overline{B}\big(\begin{smallmatrix}x&y\\
y&z\end{smallmatrix}\big)$. Avec ces notations, on voit que
\begin{equation}
d^0|_F\big(\begin{smallmatrix}\alpha &|&\beta \,\\\hline \gamma
&|&-\alpha^*\,\end{smallmatrix}\big)\big(\begin{smallmatrix}B\gamma &|&\alpha B-B\alpha^* \,\\\hline
0&|&\gamma B\,\end{smallmatrix}\big), \quad d^1|_F
\big(\begin{smallmatrix}\alpha &|&\beta \,\\\hline \gamma
&|&\alpha^*\,\end{smallmatrix}\big)\big(\begin{smallmatrix}-B\gamma&|&\alpha B-B\alpha^* \,\\\hline
0&|&\gamma B\,\end{smallmatrix}\big).
\end{equation}
Par cons\'{e}quent, $\big(\begin{smallmatrix}\alpha &|&\beta
\,\\\hline \gamma
&|&\alpha^*\,\end{smallmatrix}\big)\in\Kern(d^1|_F)$ si et
seulement si $B\gamma=0$ et  $\alpha B$ est sym\'{e}trique. De m\^{e}me
$\big(\begin{smallmatrix}\alpha &|&\beta \,\\\hline \gamma
&|&\alpha^*\,\end{smallmatrix}\big)\in\Bild(d^0|_F)$ si et
seulement si $\gamma =0$ et s'il existe $\gamma^\prime$
sym\'{e}trique avec $\alpha=B\gamma^\prime$ et $\alpha^\prime$ avec
$\beta=\alpha^\prime B-B\alpha^{\prime*}$. En particulier,
$\sigma^*C^\bull|_{F}$ se d\'{e}compose en trois parties suivant
$\alpha,\beta$ ou $\gamma$. Plus pr\'{e}cis\'{e}ment $H^1$ est la somme
directe de
\begin{enumerate}
  \item $H_{I\phantom{II}}^1=\{\alpha\ |\ \alpha B\text{ sym\'{e}trique } \}
                /\{B\gamma^\prime\ |\ \gamma^\prime \text{ sym\'{e}trique } \}$
  \item\label{item:HII}
        $H_{II\phantom{I}}^1=\{\beta\ |\ \beta \text{ anti-sym\'{e}trique }\}
                /\{\alpha^\prime B-B\alpha^{\prime*}\ |\
                \alpha^\prime\text{ quelconque }  \}$
  \item $H_{III}^1=\{\gamma\ |\ B\gamma=0 \text{ et }\gamma
                \text{ anti-sym\'{e}trique} \}$
\end{enumerate}
Commen\c{c}ons avec $(3)$. Pour $\gamma=\big(\begin{smallmatrix}0&t\\
-t&0\end{smallmatrix}\big)$ on a $B\gamma=\big(\begin{smallmatrix}-yt&xt\\
-zt&yt\end{smallmatrix}\big)$. Ainsi
\begin{equation}\label{eq:H1-III}
    H_{III}^1=\Kern(\ko_F\xra{\ (x\ y\ z)^t\ }\ko^3_F)=0
\end{equation}
Pour $(2)$ on observe que si
$\alpha^\prime    \big(\begin{smallmatrix}
        \alpha_{11}&\alpha_{12}\\
        \alpha_{21}&\alpha_{22}
      \end{smallmatrix}\big)$
on a $\alpha^\prime B-B\alpha^{\prime*}    \big(\begin{smallmatrix}
        0&v\\
        -v&0
      \end{smallmatrix}\big)$
avec $v=\alpha_{11}y+\alpha_{12}z-\alpha_{21}x-\alpha_{22}y$. Ainsi
\begin{equation}\label{eq:H1-II}
    H_{II}^1=\Kokern(\ko_F^4\xra{\ (y\ z\ -x\ y)\ }\ko_F)=\ko_{\{0\}}.
\end{equation}
Finalement pour $(1)$, on calcule avec
$\alpha    \big(\begin{smallmatrix}
        \alpha_{11}&\alpha_{12}\\
        \alpha_{21}&\alpha_{22}
      \end{smallmatrix}\big)$
et
$\gamma^\prime    \big(\begin{smallmatrix}
        \gamma_{11}&\gamma_{12}\\
        \gamma_{12}&\gamma_{22}
      \end{smallmatrix}\big)$
que $$\alpha B    \Big(\begin{smallmatrix}
        \alpha_{11}x+\alpha_{12}y&\ &\alpha_{11}y+\alpha_{12}z\\[1ex]
        \alpha_{21}x+\alpha_{22}y&\ &\alpha_{21}y+\alpha_{22}z
      \end{smallmatrix}\Big)
\quad
    B \gamma^\prime     \Big(\begin{smallmatrix}
        \gamma_{11}x+\gamma_{12}y&\ &\gamma_{12}x+\gamma_{22}y\\[1ex]
        \gamma_{11}y+\gamma_{12}z&\ &\gamma_{12}y+\gamma_{22}z
      \end{smallmatrix}\Big)
      $$
Ainsi $H_{I}^1$ est le $H^1$ du complexe
\begin{equation}\label{eq:H1-I}
    \ko_F^3\xra{\Bigg(\begin{smallmatrix}
        x&y&0\\0&x&y\\y&z&0\\0&y&z
      \end{smallmatrix}\Bigg)}\ko_F^4\xra{\ (y\  z\  -x\ y)\ }\ko_F.
\end{equation}
Maintenant, en utilisant le complexe de Koszul du point $\{0\}\in F$, on trouve
\begin{equation}\label{eq:H1-Ivrai}
    H_{I}^1=\ko_{\{xz-y^2=0\}}.
\end{equation}

On d\'{e}duit de \eqref{eq:H1-III}, \eqref{eq:H1-II} et
\eqref{eq:H1-I} que $H^1(\sigma^* C^\bull)=L\oplus M$ o\`{u} $L$ et
$M$ sont des fibr\'{e}s en droites sur $\sigma^{-1}(Z_{sing})$ et
$\sigma^{-1}(0)=G$ respectivement. Pour terminer la d\'{e}monstration,
il suffit de montrer que $M$ n'a pas de section globale ce qui
montrera que $\sigma_*M=0$.

En fait, il est cons\'{e}quence de la d\'{e}finition de $H^1_{II}$ et de
l'\'{e}quation \eqref{eq:H1-II} que
\begin{equation}
M=\Hom_{\text{anti-sym}}(\ku^*,\ku)\subset E_+\otimes\ko_G
\end{equation}
Par cons\'{e}quent, on a des isomorphismes $M\isom \Lambda^2\ku\isom
\ko_{\IP(\Lambda^2V)}(-1)$ et donc bien $\sigma_*M=0$.
\end{proof}

%%%%%%%%%%%%%%%%%%%%%%%%%%%%%%%%%%%%%%%%%%%%%%%%%%%%%%%%%%%%%%%%%%%%%%%%%%%
%%%%%%%%%%%%%%%%%%%%%%%%%%%%%%%%%%%%%%%%%%%%%%%%%%%%%%%%%%%%%%%%%%%%%%%%%%%

\section{La structure locale de l'espace de modules}
\label{sec:structure}

Soit $X$ une surface projective lisse munie d'un diviseur ample $H$.
Soit $M$ l'espace de modules des faisceaux $H$-semi-stables sur $X$
de rang et de classes de Chern fix\'{e}s. On sait que les points de $M$
sont en bijection avec les classes d'isomorphisme de faisceaux
poly-stables, i.\ e.\ somme directe de faisceaux stables. Soit $E$
un tel faisceau et soit $\Ext^p(E,E)_0$ la partie sans trace du
$\Ext^p(E,E)$. Alors on a la description suivante de la structure
locale de $M$ en $[E]$.

Soit $\widehat A$ la compl\'{e}tion de l'anneau $A=\IC[\Ext^1(E,E)]$ des
fonctions polyn\^{o}miales sur $\Ext^1(E,E)$ par rapport \`{a} l'id\'{e}al
maximal $\gothm$ des fonctions qui s'annulent dans l'origine. La
proposition suivante est bien connue:

\begin{proposition}\label{prop:kura}--- Il y a un \'{e}l\'{e}ment
$f\in\Ext^2(E,E)_0\tensor\gothm^2 \widehat A$, dite
\emph{l'application de Kuranishi}, avec les propri\'{e}t\'{e}s suivantes:
\begin{enumerate}
\item $f$ est \'{e}quivariante par rapport \`{a} l'action de
$\Aut(E)$ sur $\Ext(E,E)$ par conjugaison.
\item La partie initiale $f_2\in \Ext^2(E,E)_0\tensor\gothm^2/\gothm^3$
de $f$ est donn\'{e}e par $f_2(e)=e\cup e$ o\`{u} $e\in\Ext^1(E,E)$.
\item Soit $\gotha\subset \widehat A$ l'id\'{e}al engendr\'{e} par l'image de
l'application adjointe\\ $f:\Ext^2(E,E)_0^*\to\gothm^2 \widehat A$.
Alors il y a un isomorphisme des anneaux complets
$\widehat{\ko}_{M,[E]}\isom (\widehat A/\gotha)^{\Aut(E)}$.
\end{enumerate}
\end{proposition}

On verra $f$ ou bien comme un \'{e}l\'{e}ment de
$\Ext^2(E,E)_0\tensor\gothm^2 \widehat A$ ou bien comme une fonction
formelle $\Ext^1(E,E)\to \Ext^2(E,E)_0$ selon le contexte.
\'{E}videmment, $f$ n'est pas unique: par exemple, on pourrait y
appliquer un changement \'{e}quivariant de coordonn\'{e}es. Nous rappelons
la construction de $f$ dans l'appendice en d\'{e}rivant au m\^{e}me temps
une autre propri\'{e}t\'{e} de $f$ dont on aura besoin.

Soit $F$ un faisceau stable sur une surface $K3$ ou ab\'{e}lienne et
$V:=\Ext^1(F,F)$. La composition
\begin{equation}
\omega:V\times V\xra{\;\cup\;}\Ext^2(F,F)\xra{\;\text{trace}\;}\IC
\end{equation}
d\'{e}finit une forme symplectique non-d\'{e}g\'{e}n\'{e}r\'{e}e. Soit $E=F\oplus F$.
Alors on a $\Aut(E)=\LieGl_2$, $\Ext^1(E,E)=\Liegl_2\tensor V$ et
$\Ext^2(E,E)_0=\Liesl_2$.
 La
proposition suivante est un cas sp\'{e}cial de la proposition
\ref{prop:Ezerfaellt} de l'appendice:

\begin{proposition}\label{prop:besser}---
Il existe un \'{e}l\'{e}ment $f$ comme dans la proposition \ref{prop:kura}
tel que $f$ s'annule sur le sous-espace
\begin{equation}
\left(\begin{array}{cc}\Ext^1(F,F)&0\\0&\Ext^1(F,F)\end{array}\right)
\subset \Liegl_2\tensor\Ext^1(F,F)=\Ext^1(E,E).
\end{equation}
En fait, gr\^{a}ce \`{a} sa $\LieGl_2$--\'{e}quivariance, $f$ s'annule sur
toute l'orbite de ce sous-espace.
\end{proposition}

On suppose \`{a} partir de maintenant que $\dim(V)=4$ et on choisit une
base symplectique  $v_1,v_2,v_3,v_4$ de $V$ par rapport \`{a} laquelle
la forme $\omega$ est donn\'{e}e par la matrice
\begin{equation}
J:=\left(\begin{array}{cc|ccc}0&1&&\\-1&0&&\\\hline&&0&1\\&&-1&0
\end{array}\right).
\end{equation}
Soit $A=S^\smb (\Liegl_2\tensor V)^*$, et soit
$f:\Liesl_2^*\to\widehat A$ une application de Kuranishi
$\LieSl_2$--\'{e}quivariante pour le faisceau $E$ qui s'annule sur
$\gothd\tensor V^*$, o\`{u} $\gothd\subset\Liegl_2$ d\'{e}signe le
sous-espace de matrices diagonales. Une telle fonction existe selon
Prop.\ \ref{prop:kura} et Prop.\ \ref{prop:besser}.

On identifie $\Liegl_2\tensor V=\Liegl_2^{\oplus 4}$ en vertu de la
base $(v_1,\ldots,v_4)$ et on d\'{e}signe par $A_i:\Liegl_2^{\oplus
4}\to \Liegl_2$, $i=1,..,4$, la $i$-i\`{e}me projection. On consid\`{e}re
$A_1,\ldots,A_4$ comme des coordonn\'{e}es \`{a} valeurs matricielles.
Alors, comme la partie quadratique de $f$ est donn\'{e}e par le produit
de Yoneda on trouve que $f$ s'exprime en termes de ces $A_i$ comme
suit:
\begin{eqnarray*}
f(A_1,A_2,A_3,A_4)&=&[A_1,A_2]+[A_3,A_4]\\
&&+\text{ \tdph }\in \Liesl_2,\text{ et}
\end{eqnarray*}
\begin{eqnarray*}
f(A_1,A_2,A_3,A_4)&=&0,\text{ si les $A_i$ sont des matrices
diagonales.}
\end{eqnarray*}
Soit $L\subset A$ l'id\'{e}al de l'adh\'{e}rence du $\LieSl_2$--orbite de
$\gothd^{\oplus 4}$ dans $\Liegl_2^{\oplus 4}$. On a donc $\gotha
\subset L\widehat A$, o\`{u} $\gotha\subset \widehat A$ est l'id\'{e}al de
la proposition \ref{prop:kura}.

Selon le premier et le deuxi\`{e}me th\'{e}or\`{e}me de la th\'{e}orie des
invariants \cite{Weyl} pour le groupe $\LieSl_2$ le sous-anneau des
invariants $A^{\LieSl_2}$ est engendr\'{e} par les fonctions suivantes:
($i,j=1,\ldots,4$)
\begin{enumerate}
\item $X_i=\trace(A_i)$,%\quad $i=1,\ldots,4$,
\item $Y_{ij}=\trace(A_i'A_j')$, %\quad $i,j=1,\ldots,4$,
o\`{u} $A_i'=A_i-\tfrac12\, X_i\,\id$,
\item $T_{i}=(-1)^i\trace(A_1',\ldots, \widehat{A'_i},\ldots,A'_4)$, %$i=1,\ldots,4$,
o\`{u}~ $\widehat{ } $~ signifie qu'on enl\`{e}ve le facteur en question.
\end{enumerate}
Si $Y=(Y_{ij})$ et $T=(T_i)$, alors les relations fondamentales
entre elles sont:
\begin{equation}
Y=Y^t,\quad \det(Y)=0,\quad YT=0,\quad TT^t=-2\text{adj}(Y),
\end{equation}
o\`{u} $\text{adj}(Y)$ d\'{e}signe la matrice adjointe de $Y$.
% soumis seulement aux r\'{e}lations $Y_{ij}=Y_{ji}$, $\det(Y)=0$ (o\`{u} $Y$
% designe la matrice sym\'{e}trique $(Y_{ij})$) et
% $T_{i_1i_2i_3}T_{j_1j_2j_3}=\det( (Y_{i_aj_b})_{ab})$.

D'autre part, l'id\'{e}al $\gotha_0=\gotha\cap\widehat A^{\LieSl_2}$
contient la fonction
\begin{equation}\label{eq:Tijk}
T'_{1}:=-\frac12\trace(A_2' f(A_1,\ldots,A_4))=T_{1}+ \text{\tdph}
\end{equation}
(et d'une fa\c{c}on similaire des fonctions $T_{i}'=T_{i}+\ldots$ pour
les autres indices), et aussi
\begin{equation}%
\trace(A_i'A_j'f(A_1,\ldots,A_4))= (YJY)_{ij}+\text{\tdph}.
\end{equation}
Gr\^{a}ce \`{a} \eqref{eq:Tijk}, on peut remplacer les $T_{i}$ par les
$T_{i}'$ comme coordonn\'{e}es. Ensuite, en calculant modulo $\gotha_0$,
on peut \'{e}liminer les $T_{i}'$ compl\`{e}tement. On en d\'{e}duit que
\begin{equation}%
(\widehat A/\gotha)^{\LieSl_2}=\widehat A^{\LieSl_2}/\gotha_0\isom
\IC[[X_1\ldots,X_4,Y_{11},Y_{12},\ldots,Y_{44}]]/I
\end{equation}
pour un certain id\'{e}al $I$ qui contient des fonctions
\begin{equation}%
f_{ab}(X_i,Y_{jk}):= (YJY)_{ab} +\text{\tdph},\quad 1\leq a\leq b\leq 4.
\end{equation}
Nous laissons au lecteur le soin de v\'{e}rifier que $L\cap
A^{\LieSl_2}$ est engendr\'{e} par tous les $2\times 2$--mineurs
$Y_{ab}Y_{cd}-Y_{ad}Y_{bc}$ de $Y$.

Pour simplifier les notations dans la discussion suivante on
identifie l'espace des matrices sym\'{e}trique $\Sym_4$ avec l'alg\`{e}bre
de Lie symplectique $\Liesp_4$ en vertu de l'application $Y\mapsto
B:=YJ$. Par abus de notation on garde la notation $I$. Rappelons que
$I_0\subset L_0\subset\IC[\Liesp_4]$ d\'{e}signe l'id\'{e}aux engendr\'{e} par
les coefficients des matrices $B^2$ et $\Lambda^2B$, respectivement
(Cf.\ Sec.\ \ref{section:resolution}). R\'{e}sumons donc:

\begin{proposition}\label{prop:resume}--- Soit $\widehat R$ la compl\'{e}tion de
$R= \IC[\IC^4\times \Liesp_4]$ dans l'origine. Il y a des id\'{e}aux
$I'\subset I\subset L_0\widehat R\subset\widehat R$ tels que
\begin{enumerate}
\item[(i)] il y a un isomorphisme
$\Phi:\widehat\ko_{M,[E]}\isom \widehat A^{\LieSl_2}/\gotha_0\isom\widehat R/I$,
\item[(ii)] $I'$ est engendr\'{e} par des fonctions $f_{\alpha}(X,B)$,
$\alpha=1,\ldots,6$, dont les termes initiaux quadratiques engendrent l'id\'{e}al
$I_0 R$. Ici  $X\in\IC^4$ et $B\in\Liesp_4$.
\item[(iii)] Sous l'isomorphisme $\Phi$ l'id\'{e}al $L_0$ correspond au lieu
des faisceaux strictement semi-stables (= semi-stables mais non stables), et
donc au lieu singulier du germe formel.
\end{enumerate}
\end{proposition}

\begin{proof} (i) est cons\'{e}quence de la discussion pr\'{e}c\'{e}dente.
Pour voir (ii), on remarque que $I$ contient les fonctions $\tilde
f_{ac}=\sum_b f_{ab}J_{bc}$ dont les termes quadratiques sont
donn\'{e}es par $(YJYJ)_{ac}=(B^2)_{ac}$. L'id\'{e}al $I_0$ est engendr\'{e} par
6 quadriques ; ainsi on peut choisir les $f_\alpha$ parmi les
$\tilde f_{ac}$. Pour (iii) et l'inclusion $I\subset L_0\widehat R$,
il suffit de v\'{e}rifier que $L_0 R=L\cap A^{\LieSl_2}$ qu'on voit par
un calcul direct.
\end{proof}

Soit $h\in\widehat R$ une fonction non nulle et
$h=h_0+h_1+h_2+\ldots$ sa d\'{e}composition en termes homog\`{e}nes. L'ordre
de $h$ est l'index minimal $\nu$ tel que $h_\nu\neq 0$. Le terme
initial de $h$ est $\ini(h):=h_\nu\in R$. Les termes initiaux
$\ini(h)$ de tous les \'{e}l\'{e}ments $h$ d'un id\'{e}al $\gothb\subset\widehat
R$ engendrent son id\'{e}al initial $\ini(\gothb)\subset R$. Il a y un
isomorphisme canonique d'anneaux gradu\'{e}s associ\'{e}s $\gr(\widehat
R/\gothb)=R/\ini(\gothb)$.

% Soit $I_0\subset R$ l'id\'{e}al engendr\'{e} par les coefficients de la
% matrice $B^2$. Un calcul rapides montre que les 16 \'{e}quations
% correspondants se r\'{e}duisent effectivement \`{a} six quadriques, et la
% v\'{e}rification que $I_0$ est un id\'{e}al premier de hauteur 4 se fait le
% plus simplement  par ordinateur, par exemple avec \textsc{Singular}.

\begin{lemma}\label{lem:Flachheit}--- Avec les notations de proposition
\ref{prop:resume} on a
\begin{equation}
I'=I\quad\quad\text{et}\quad\quad\ini(I)=I_0 R.
\end{equation}
\end{lemma}

\begin{proof} Par d\'{e}finition de $I'$,  on a les inclusions
$I_0 R\subset \ini(I')\subset \ini(I)$. Ils induisent de surjections
$\gr(\widehat R/I)\to \gr(\widehat R/I')\to R/I_0R$ et impliquent
ainsi  des in\'{e}galit\'{e}s
$$10=\dim_{[E]}M=\dim \widehat \ko_{M,[E]}=\dim(\widehat R/I)\geq
\dim(\widehat R/I')\geq \dim R/I_0.$$ Comme $\dim R/I_0R=10$ on doit
avoir \'{e}galit\'{e} partout, et comme $I_0$ est un id\'{e}al premier il
s'ensuit que $I_0R=\ini(I')=\ini(I)$. Si donc les termes initiaux
des $f_\alpha$ engendrent l'id\'{e}al initial de $I$, les $f_\alpha$
elles-m\^{e}mes engendrent $I$, c'est-\`{a}-dire $I'=I$.
\end{proof}

\begin{theorem}\label{thm:isomorphisme}--- Soit $X$ une surface K3
ou ab\'{e}lienne et $F$ un faisceau stable par rapport \`{a} une
polarisation $H$ tel que $\dim Ext^1(F,F)=4$. Soit $E=F\oplus F$ et
$M$ l'espace de modules des faisceaux semi-stables ayant vecteur de
Mukai $v(E)$. Alors il y a un isomorphisme des germes d' espaces
analytiques
\begin{equation}
(M,[E])\isom(\Spec_{{\rm an}}R/I_0R,0)\isom (\IC^4\times Z,0).
\end{equation}
\end{theorem}

Pour d\'{e}montrer ce th\'{e}or\`{e}me, il suffira, gr\^{a}ce au
th\'{e}or\`{e}me d'Artin (voir \cite{Artin} Cor.\ 1.6), de trouver un
isomorphisme $\widehat\ko_{M,[E]}\isom (R/I_0)^\wedge$ d'anneaux complets.
En vue de la proposition \ref{prop:resume} on est ramen\'{e} \`{a} prouver que
\begin{equation}\label{eq:dasbleibtzuzeigen}
\widehat R/I=\widehat R/I_0\widehat R.
\end{equation}
On va prouver cet isomorphisme en montrant que la d\'{e}formation de $\widehat R/I$
vers son c\^{o}ne normal est essentiellement une d\'{e}formation triviale.

%%%%%%%%%%%%%%%%%%%%%%%%%%%%%%%%%%%%%%%%%%%%%%%%%%%%%%%%%%%%%%%%%%%%%%%%%%%%%
%%%%%%%%%%%%%%%%%%%%%%%%%%%%%%%%%%%%%%%%%%%%%%%%%%%%%%%%%%%%%%%%%%%%%%%%%%%%%

\section{La d\'{e}formation vers le c\^{o}ne normal}
\label{sec:conenormal}

On garde les notations de la section pr\'{e}c\'{e}dente.

Soit $\gothm\subset R$ l'id\'{e}al maximal des fonctions qui s'annulent
\`{a} l'origine et soit $R[t]^\wedge$ la compl\'{e}tion par rapport \`{a} la
topologie $\gothm R[t]$--adique. Comme ci-dessus on d\'{e}signe par
\begin{equation}%
f_\alpha=f_{\alpha,2}+f_{\alpha,3}+f_{\alpha,4}+\ldots,
\quad \alpha=1,\ldots,6,
\end{equation}
la d\'{e}composition en termes homog\`{e}nes et on pose
\begin{equation}%
F_\alpha(t):=f_{\alpha,2}+tf_{\alpha,3}+t^2f_{\alpha,4}+\ldots,\quad
\alpha=1,\ldots,6.
\end{equation}
Le quotient $S=R[t]^\wedge/\kit$ avec $\kit:=(F_1,\ldots,F_6)$ satisfait \`{a}
\begin{equation}%
S_\lambda:=S/(t-\lambda)S=\left\{
\begin{array}{ll}%
\widehat R/I_0\widehat R,& \text{ si }\lambda=0, et\\
\widehat R/I,&\text{ si } \lambda\neq 0.
\end{array}\right.
\end{equation}
Gr\^{a}ce au fait que $I_0R=\ini(I)$ (Lemme \ref{lem:Flachheit}), $S$
est plat sur $\IC[t]$.

A partir de maintenant on \'{e}crira $x_1,\ldots,x_{14}$ pour les
coordonn\'{e}es $X_i$ et $B_{ab}$ de $\IC^4\times\Liesp_4$.

On dira qu'un \'{e}l\'{e}ment de $R[t]^\wedge$ est homog\`{e}ne de poids $k$, si
son coefficient devant $t^m$ est homog\`{e}ne de degr\'{e} $k+m$. Par
exemple, $\partial_{x_i}F_\alpha$, $F_\alpha$ et $\partial_t
F_\alpha$ sont homog\`{e}nes de poids $1$, $2$ respectivement $3$.
\'{E}videmment, un \'{e}l\'{e}ment arbitraire $g\in R[t]^\wedge$ peut s'\'{e}crire
comme s\'{e}rie convergente $g=\sum_{k\in\IZ} g^{(k)}$ avec des \'{e}l\'{e}ments
homog\`{e}nes  uniques $g^{(k)}\in R[t]^\wedge$ de poids $k$.

\begin{proposition}\label{prop:ilyaPhiEtH}--- Il y a des fonctions
homog\`{e}nes $\Phi_i(x,t)\in R[t]^\wedge$ de poids $2$ pour
$i=1,\ldots,14$, et des fonctions $h_{\alpha\beta}(x,t)\in
R[t]^\wedge$ de poids $1$ pour $\alpha,\beta=1,\ldots,6$, telles que
\begin{equation}\label{eq:ft}
\sum_i \Phi_i(x,t) \partial_{x_i}F_\alpha(x,t)=\partial_t
F_\alpha(x,t) +\sum_\beta h_{\alpha\beta}(x,t)F_\beta(x,t).
\end{equation}
\end{proposition}

\begin{proof} On consid\`{e}re le module $T^1_{S/\IC[t]}$ d\'{e}fini
par la suite exacte
\begin{equation}%
S\langle\partial_{x_1},\ldots,\partial_{x_{14}}\rangle
\xra{\;\Theta\;}
Hom_{S}(\kit/\kit^2,S)\lra T^1_{S/\IC[t]}\lra 0
\end{equation}
et de m\^{e}me le module $T^1_{S_0}$ d\'{e}finit par
\begin{equation}%
S_0\langle\partial_{x_1},\ldots,\partial_{x_{14}}\rangle\xra{\;\Theta\;}
Hom_{\ko_{S_0}}(I_0/I_0^2,\ko_{S_0})\lra T^1_{S_0}\lra 0,
\end{equation}
o\`{u} $\Theta$ signifie l'action d'un champ vectoriel sur les fonctions
de $\kit$ et de $I_0$ respectivement. Consid\'{e}rons maintenant
l'\'{e}l\'{e}ment
\begin{equation}
\partial_t\in Hom(\kit/\kit^2,\ko_S),\quad\partial_t(F_\alpha)
=f_{\alpha,3}+2t f_{\alpha,4}+3t^2f_{\alpha,5}+\ldots.
\end{equation}
Trouver des fonctions $\Phi_i$ et $h_{\alpha\beta}$ satisfaisant
l'\'{e}quation \eqref{eq:ft} \'{e}quivaut \`{a} montrer que $\partial_t$ est
contenu dans l'image de l'application $\Theta$ ou bien que la classe
de $\partial_t$ dans $T^1_{S/\IC[t]}$ est nulle.

On utilise les faits suivants: Par le crit\`{e}re de Jacobi le module
$T^1_{S/\IC[t]}$ est support\'{e} sur le lieu singulier relative de $S$.
Par les r\'{e}sultats de O'Grady \cite{OG1}, les singularit\'{e}s de
$S_\lambda$ sont de type $A_1$ en dehors de l'origine. En
particulier, $T^1_{S/\IC[t]}$ est un module inversible en dehors de
l'origine et donc est annul\'{e} par l'id\'{e}al $L_0$ de lieu singulier.
D'autre part, comme $I\subset L_0$ et comme $L_0$ est engendr\'{e} par
des \'{e}l\'{e}ments homog\`{e}nes, tous les composants homog\`{e}nes des $f_\alpha$
sont contenus dans $L_0$. Par cons\'{e}quent, $\partial_t$ est annul\'{e}
par $L_0$ et d\'{e}finit donc une section de $T^1_{S/\IC[t]}$ qui, pour
$t=\lambda$ fix\'{e}, est support\'{e}e \`{a} l'origine. Il suffit donc de
montrer que $T^1_{S/\IC[t]}$ est un module pur, c'est-\`{a}-dire sans
sous-module de dimension plus petite. La multiplication par le
param\`{e}tre $t$ fournit la suite exacte
\begin{equation}%
T^1_{S/\IC[t]}\xra{\;t\;}T^1_{S/\IC[t]}\lra T^1_{S_0}.
\end{equation}
On d\'{e}duit du th\'{e}or\`{e}me \ref{thm:purete} que $T^1_{S_0}$ est pur de
dimension $8$ et g\'{e}n\'{e}riquement de rang 1. Il s'ensuit que
$T^1_{S/\IC[t]}$ est pur de dimension $9$ et donc que la classe de
$\partial_t$ dans $T^1_{S_0}$ est nulle.

Ainsi il y a des fonctions $\Phi_i$, $h_{\alpha\beta}\in
R[t]^\wedge$ telles qu'on ait \eqref{eq:ft}. Il reste \`{a} montrer que
l'on peut choisir ces fonctions comme homog\`{e}nes. Mais comme
$F_\alpha$ et ses d\'{e}riv\'{e}es sont homog\`{e}nes, les parties homog\`{e}nes de
degr\'{e} $2$ dans $\Phi_i$ et de degr\'{e} $1$ dans $h_{\alpha\beta}$ sont
aussi des solutions de \eqref{eq:ft}. Donc on peut supposer que
$\Phi_i$ et $h_{\alpha\beta}$ sont homog\`{e}nes.
\end{proof}

\begin{proposition}--- Il y a des fonctions $\Psi_i\in R[t]^\wedge$ homog\`{e}nes
de poids $1$ pour $i=1,\ldots,14$, et $M_{\alpha\beta}\in
R[t]^\wedge$ homog\`{e}nes de poids $0$ pour $\alpha,\beta=1,\ldots,6$,
telles que
\begin{enumerate}
\item $M_{\alpha,\beta}(x,0)=\delta_{\alpha\beta}$,
\item $\Psi_i(x,0)=x_i$,
\item $F_\alpha(\Psi(x,t),t)=\sum_\beta M_{\alpha\beta}(\Psi(x,t),t)f_{\beta,2}(x)$.
\end{enumerate}
\end{proposition}

\begin{proof} Soient $\Phi_i$ et $h_{\alpha\beta}$ les fonctions de la proposition
\ref{prop:ilyaPhiEtH}. On d\'{e}finit un changement de coordonn\'{e}es
formel $\Psi(x,t)$ et une matrice inversible $M(x,t)$ comme
solutions du syst\`{e}me des \'{e}quations diff\'{e}rentielles suivantes:
\begin{equation}\label{eq:Psit}
\left\{\begin{array}{l}
\partial_t\Psi_i(x,t)=-\Phi_i(\Psi(x,t),t),\\
\Psi_i(x,0)=x_i,
\end{array}
\right.
\end{equation}
et
\begin{equation}\label{eq:Mt}
\left\{\begin{array}{l}
\partial_t M_{\alpha\beta}(x,t)=\sum_i\Phi_i(x,t)\partial_iM_{\alpha\beta}(x,t)
 -\sum_{\beta}h_{\alpha\gamma}(x,t)M_{\gamma\beta}(x,t)\\
M_{\alpha\beta}(x,0)=\delta_{\alpha\beta}
\end{array}\right.
\end{equation}
Remarquons que -- gr\^{a}ce au fait que $\Phi_i$ et $h_{\alpha\beta}$
sont homog\`{e}nes et commencent avec des termes quadratique
respectivement lin\'{e}aires en les variables $x_i$ --  ses \'{e}quations se
r\'{e}solvent ais\'{e}ment par r\'{e}currence par rapport au degr\'{e} en les
variables $x_i$.

En utilisant les relations \eqref{eq:ft}, \eqref{eq:Psit} et
\eqref{eq:Mt} on v\'{e}rifie qui les fonctions
\begin{equation}%
g_\alpha(x,t):=F_\alpha(\Psi(x,t),t)-\sum_\beta
M_{\alpha\beta}(\Psi(x,t),t) f_{\beta,2}(x)
\end{equation}
satisfont \`{a} l'\'{e}quation diff\'{e}rentielle
\begin{equation}%
\partial_tg_\alpha(x,t)=-\sum_\beta h_{\alpha\beta}(\Psi(x,t),t)g^\beta(x,t).
\end{equation}
Comme $g_\alpha(x,0)=0$, on d\'{e}duit que $g_\alpha(x,t)=0$.
\end{proof}

Avec les notations de la proposition pr\'{e}c\'{e}dente on pose
\begin{equation}
\Psi(x):=(\Psi_i(x,1))\in \Aut(\widehat R)
\end{equation}
et
\begin{equation}
M(x):=(M_{\alpha\beta}(x,1))\in \LieGl_6(\widehat R).
\end{equation}
En \'{e}valuant les identit\'{e}s de la proposition en $t=1$ on obtient:
\begin{equation}
f_\alpha(\Psi(x))=\sum_\beta M_{\alpha\beta}(x) f_{\beta,2}(x).
\end{equation}
Donc l'automorphisme $\Psi:\widehat R\to\widehat R$ envoie $I$ sur
$I_0\widehat R$ et induit un isomorphisme
\begin{equation}
\Psi:\widehat R/I\lra \widehat R/I_0\widehat R.
\end{equation}
ce qui termine la d\'{e}monstration de \eqref{eq:dasbleibtzuzeigen} et
donc du th\'{e}or\`{e}me \ref{thm:isomorphisme}.
\bigskip

\noindent {\sl Fin de la d\'{e}monstration du th\'{e}or\`{e}me} \ref{th:main}:
Rappelons que le vecteur de Mukai d'un faisceau coh\'{e}rent $E$ sur une
surface $K3$ projective ou ab\'{e}lienne est donn\'{e} par
$v(E)=\ch(E)\sqrt{\text{td}(X)}$ (voir \cite{HL}). Alors se donner
un vecteur de Mukai \'{e}quivaut \`{a} fixer le rang et les classes de Chern
d'un faisceau, et $\dim \Ext^1(E,E)= \langle v(E),v(E)\rangle +2
\dim\End(E)$. Les hypoth\`{e}ses du th\'{e}or\`{e}me entra\^{\i}nent qu'il y a
seulement trois types de faisceaux poly-stable dans $M_{2v}$:
\begin{enumerate}
\item[(1)] Les faisceaux stables $E$ avec $v(E)=v$. Par les r\'{e}sultats de Mukai,
$M_{2v}$ est lisse dans $[E]$ de dimension \'{e}gale \`{a}
$\dim\Ext^1(E,E)=4\cdot 2+2=10$.
\item[(2)] Les faisceaux $E=F_1\oplus F_2$ avec des faisceau stable non-isomorphes
$F_1$ et $F_2$ o\`{u} $v(F_1)=v(F_2)=v$. Chaque composante $F_i$ varie dans
l'espace de module $M_v$ de dimension $4$, et donc les points dans $M_{2v}$
de ce type forme une vari\'{e}t\'{e} $S^2M_v\setminus\Delta_{M_v}$ de dimension
$8$. Ils sont les points singuliers g\'{e}n\'{e}riques de $M_{2v}$.
\item[(3)] Les faisceaux $E=F^{\oplus2}$ avec un faisceau stable $F$ de vecteur de
Mukai $v(F)=v$. Gr\^{a}ce \`{a} l'hypoth\`{e}se $\langle v,v\rangle=2$, l'espace
vectoriel symplectique $V=\Ext^1(F,F)$ a la dimension $4$.
\end{enumerate}
De cette fa\c{c}on, O'Grady \cite{OG1} a obtenu une stratification
$M_{2v}\supset S^2M_v\supset\Delta_{M_{v}}$ de l'espace de modules
singulier. Il a aussi montr\'{e} que la singularit\'{e} d'un point de
$S^2M_v\setminus\Delta_{M_v}$ est de type $A_1$ \`{a} travers de
$S^2M_v$. Donc l'\'{e}clatement de $M_{2v}\setminus\Delta_{M_v}$ le long
$S^2M\setminus \Delta_{M_v}$ est une r\'{e}solution symplectique.

Les faisceaux de type (3) correspondants aux points de
$\Delta_{M_v}$ satisfont aux hypoth\`{e}ses des th\'{e}or\`{e}mes
\ref{thm:isomorphisme} et \ref{thm:laresolution}. On d\'{e}duit que
l'\'{e}clatement de $M_{2v}$ le long de $S^2M_v$ est une r\'{e}solution
semi-petite. \qed

%%%%%%%%%%%%%%%%%%%%%%%%%%%%%%%%%%%%%%%%%%%%%%%%%%%%%%%%%%%%%%%%%%%%%%%%%%%%%%%%
%%%%%%%%%%%%%%%%%%%%%%%%%%%%%%%%%%%%%%%%%%%%%%%%%%%%%%%%%%%%%%%%%%%%%%%%%%%%%%%%

\setcounter{section}{0}
\renewcommand{\thesection}{\Alph{section}}
\section{Appendice: L'application de Kuranishi}

Soit $E$ un faisceau polystable sur une surface projective lisse.
Les d\'{e}formations infinit\'{e}simales de $E$ se d\'{e}crivent en termes des
solutions de l'\'{e}quation de Maurer-Cartan comme suit:

Soit $E=\bigoplus_iE_i^{\oplus n_i}$ la d\'{e}composition de $E$ comme
somme directe de faisceaux stables. En choisissant des r\'{e}solutions
$E_i\to I_{(i)}$ et en formant $I=\bigoplus_i I_{(i)}^{\oplus n_i}$
on obtient une r\'{e}solution $E\to I$ qui est injective  et
\'{e}quivariante pour l'action canonique du groupe
$\Aut(E)=\prod_i\LieGl_{n_i}$.

Le diff\'{e}rentiel $d_I$ de $I$ est un \'{e}l\'{e}ment de degr\'{e} 1 dans
$\gothg:=\End(I,I)$ et d\'{e}finit un diff\'{e}rentiel \'{e}quivariant
$d=\ad(d_I)$ sur $\gothg$. Alors $H(\gothg,d)=\Ext(E,E)$. Soient
$B\gothg\subset Z\gothg\subset\gothg$ les sous-complexes des cobords
et des cocycles. On choisit, une fois pour tout, des homomorphismes
\'{e}quivariants $s:H\gothg\to Z\gothg$ et $t:B\gothg\to\gothg$ de degr\'{e}
0 et $-1$ respectivement qui scindent les surjections naturelles
$\overline{\phantom{m}}:Z\gothg\to H\gothg$ et $d:\gothg\to
B\gothg$.

Soit $(A',\gothm')$ un anneau local artinien avec $A'/\gothm'=\IC$
et soit $d'\in \gothg^1\tensor A'$ tel que
\begin{equation}
d'\tensor A'/\gothm'=d_I\quad\text{et}\quad d'{}^2=0. \end{equation}
On sait que $(I\tensor A',d')$ est un complexe acyclique et que
$E':=H^0(I\tensor A',d')$ est une d\'{e}formation plate de $E$ sur $A'$.
R\'{e}ciproquement, toutes les d\'{e}formations s'obtiennent ainsi.
Supposons en plus que $A''\to A'$ est un \'{e}pimorphisme des anneaux
artiniens dont le noyau $J\subset A''$ est annul\'{e} par l'id\'{e}al
maximal $\gothm''\subset A''$. Pour \'{e}tendre $d'$ \`{a} $A''$ il faut
choisir un \'{e}l\'{e}ment $d''\in\gothg\tensor A''$ qui se restreint \`{a}
$d'$. Gr\^{a}ce \`{a} la relation $d'{}^2=0$ on a $d''{}^2\in\gothg^2\tensor
J$. Cet \'{e}l\'{e}ment est en fait un cocycle et sa classe
$\gotho(d',A'')\in \Ext^2(E,E)\tensor J$ ne d\'{e}pend pas du choix de
$d''$. Il est aussi bien connu (voir \cite{Artamkin,Artamkin2}) que
la trace de cette classe d'obstruction est nulle. Donc le plus grand
quotient $A''/\gotha$ tel qu'une extension de $E'$ \`{a} $A''$ existe
est donn\'{e} par l'id\'{e}al $\gotha$ engendr\'{e} par l'image de l'application
ajointe $\gotho(d',A''):\Ext^2(E,E)_0^*\to J\subset A''$.

On est amen\'{e} \`{a} l'algorithme suivant pour le calcul de la d\'{e}formation
verselle de $E$: Soit $U=\Ext^1(E,E)^*$, $A=S^\smb
U=\IC[\Ext^1(E,E)]$, $\widehat A$ la compl\'{e}tion de $A$ \`{a} l'origine
et $\gothm\subset \widehat A$ son id\'{e}al maximal. On d\'{e}signe par
$\gamma_1=(s\tensor 1)(\id_U) \in Z^1(\gothg\tensor \widehat A)$ le
cocycle tautologique.

\begin{proposition}\label{prop:kuraconst}---
Il y a des \'{e}l\'{e}ments \'{e}quivariants
$$\gamma_n\in \gothg^1\tensor S^nU\quad\text{et}\quad
f_n\in \Ext^2(E,E)_0\tensor S^nU,\quad\text{pour }n\geq 2,$$ avec
les propri\'{e}t\'{e}s suivantes: Si $\gamma=\gamma_1+\gamma_2+\ldots$ et
$f=f_2+f_3+\ldots$, et si on d\'{e}signe par $\gotha\subset
\gothm^2\subset \widehat A$ l'id\'{e}al engendr\'{e} par l'image de
$f:\Ext^2(E,E)_0^*\to \widehat A$, alors
$$(d_I+\gamma)^2-s(f)\in\gothg^2\tensor\gotha\gothm.$$
En particulier, $d_I+\gamma$ d\'{e}finit une d\'{e}formation plate de $E$
sur $\widehat A/\gotha$.
\end{proposition}

\begin{proof} On va construire $f_n$ et $\gamma_n$ par r\'{e}currence.
On suppose que $n\geq 2$ et qu'on a d\'{e}j\`{a}
trouv\'{e} des \'{e}l\'{e}ments $\gamma_2,\ldots,\gamma_{n-1}$ et
$f_2,\ldots,f_{n-1}$ tels que
\begin{equation}
x:=(d+\gamma_1+\ldots+\gamma_{n-1})^2-s(f_2+\ldots+f_{n-1})
\in\gothg^2\tensor (\gotha\gothm+\gothm^n).
\end{equation}
Comme $\gamma_1$ est un cocyle, cette hypoth\`{e}se est certainement
v\'{e}rifi\'{e}e pour $n=2$. Remarquons que $\gotha\gothm+\gothm^{n+1}$ est
bien d\'{e}finie par $f_2,\ldots,f_{n-1}$. Or,
\begin{equation}
d(x)\in\gothg^3\tensor(\gotha\gothm+\gothm^{n+1})
\end{equation}
puisque
\begin{eqnarray*}
0&=&[d+\gamma_1+\ldots+\gamma_{n-1},(d+\gamma_1+\ldots+\gamma_{n-1})^2]\\
&=&d(x)+[\gamma_1+\ldots+\gamma_{n-1},s(f_2+\ldots+f_{n-1})+x].
\end{eqnarray*}
Donc la classe
$$[x]\in \gothg^2\tensor\frac{\gotha\gothm+\gothm^{n}}{\gotha\gothm+\gothm^{n+1}}
$$
de $x$ est un cocycle. Sa classe de cohomologie $\overline{[x]}$ est
l'obstruction pour l'extension de  $\widehat A/\gotha+\gothm^{n}$ \`{a}
$\widehat A/\gothm(\gotha+\gothm^{n})$ et donc sans trace. Selon le
th\'{e}or\`{e}me de Weyl il y a un scindage \'{e}quivariant $q$ de la projection
de $\Aut(E)$--modules
$$S^n U =\frac{\gothm^{n}}{\gothm^{n+1}}\lra \frac{\gotha\gothm+\gothm^{n}}{\gotha\gothm+\gothm^{n+1}}.
$$
Soit $z=q([x])\in Z^2\gothg\tensor S^nU$. On pose
$f_n=\overline{z}$, $\gamma_n:=t(s(f_n)-z)\in\gothg^1\tensor S^nU$
et $x':=x-z+[\gamma_1+\ldots+\gamma_n,\gamma_n]$. Par construction,
$f_n$ et $\gamma_n$ sont \'{e}quivariants, et
$$(d+\gamma_1+\ldots+\gamma_n)^2-s(f_2+\ldots+f_n)x'\in \gothg^2\tensor(\gotha\gothm+\gothm^{n+1}).$$
\end{proof}

Soit maintenant $X$ une surface $K3$. Soient $F$ un faisceau stable
et $E=F^{\oplus r}$. Alors $\Aut(E)\isom\LieGl_r$ et la r\'{e}solution
$E\to I_E$ a la forme $I_E\isom\IC^r\tensor I_F$ pour une r\'{e}solution
$F\to I_F$. Aussi, $\End(I_E,I_E)\isom\Liegl_r\tensor \gothg$ avec
$\gothg=\End(I_F)$.

On d\'{e}note par $F_a$ la $a$-i\`{e}me composante de $E$. Il y a une
d\'{e}composition canonique $U=\bigoplus_{ab}U_{ab}$ avec
$U_{ab}=\Ext^1(F_b,F_a)^*$. Soit $\gothq\subset (S^\smb U)^\wedge$
l'id\'{e}al engendr\'{e} par les sous-espace $U_{ab}\subset U$, $a\neq b$,
et soit $U':=U_{11}\oplus\ldots\oplus U_{rr}$, tel que $(S^\smb
U)^\wedge/\gothq=(S^\smb U')^\wedge$.

\begin{proposition}\label{prop:Ezerfaellt}---
Soient $\gamma$ et $f$ des \'{e}l\'{e}ments construits
pour le faisceau $E$ par la m\'{e}thode de la d\'{e}monstration pr\'{e}c\'{e}dente.
Alors\begin{small}
$$f\left(\begin{smallmatrix}
e_1&0&&0\\
0&e_2&&\\
&&\ddots&\\
0&&&e_r
\end{smallmatrix}\right)
\equiv
\left(\begin{smallmatrix}
0&&&\\
&0&&\\
&&\ddots&\\
&&&0
\end{smallmatrix}\right)
\text{ et }
\gamma
\left(\begin{smallmatrix}
e_1&0&&0\\
0&e_2&&\\
&&\ddots&\\
0&&&e_r
\end{smallmatrix}\right)
\equiv
\left(\begin{smallmatrix}
\tilde\gamma(e_1)&0&&0\\
0&\tilde\gamma(e_2)&&\\
&&\ddots&\\
0&&&\tilde\gamma(e_r)
\end{smallmatrix}\right)
$$\end{small}pour $e_a\in \Ext^1(F_a,F_a)$, $a=1,\ldots r$, o\`{u} $\tilde\gamma\in\gothg^1
\tensor (S^\smb \Ext^1(F,F)^*)^\wedge$ satisfait \`{a} $(d+\tilde\gamma)^2=0$.
\end{proposition}

\begin{proof} La proposition affirme que $f$ s'annule modulo $\gothq$ et que
-- encore modulo $\gothq$ -- $\gamma$ prend une forme diagonale avec
coefficients $\gamma_{aa}\in (S^\smb U_{aa})^\wedge$. Que ces
coefficients sont donn\'{e}es par la m\^{e}me fonction $\tilde\gamma$ est
d'ailleurs une cons\'{e}quence imm\'{e}diate de l'\'{e}quivariance de $\gamma$
sous l'action du groupe sym\'{e}trique.

On va montrer par induction que l'affirmation est vraie modulo
$\gothm^n$ pour tous  $n\geq 2$. Le cas $n=2$ est trivial gr\^{a}ce \`{a} la
nature tautologique du cycle $\gamma_1$. Supposons alors que $n>2$
et que le r\'{e}sultat est d\'{e}j\`{a} montr\'{e} modulo $\gothm^n$. Donc par
hypoth\`{e}se d'induction $\gotha\subset\gothq+\gothm^n$. Cela implique
que la fl\`{e}che $i$ du diagramme canonique
\begin{equation}\label{eq:vergleich}
\begin{array}{ccc}
S^nU&\lra&\frac{\gothm^n}{\gothm^{n+1}+\gotha\gothm}\\
\Big\downarrow&&\Big\downarrow\\
S^nU'&\xra{\;\;i\;\;}&
\frac{\gothm^n}{\gothm^{n+1}+\gotha\gothm+\gothq\cap\gothm^{n}}
\end{array}
\end{equation}
est un isomorphisme.

On garde les notations de la d\'{e}monstrations pr\'{e}c\'{e}dente. Parce que
l'\'{e}l\'{e}ment $x$ s'exprime compl\`{e}tement en termes de $f$ et de $\gamma$
il suit que la classe de $x$ modulo $\gothq$ a la m\^{e}me structure
diagonale comme $\gamma$. Il suit de la commutativit\'{e} du diagramme
\eqref{eq:vergleich} et du fait que $i$ est un isomorphisme que la
projection de $z\in S^nU$ dans $S^nU'$ a aussi une telle structure,
c'est-\`{a}-dire est contenue dans le sous-espace
$S^nU_{11}\oplus\ldots\oplus S^nU_{rr}$. Cette propri\'{e}t\'{e} est h\'{e}rit\'{e}e
par $f_n =\overline{z}$ et $\gamma_n=s(f_n)-z_n$. En particulier,
les coefficients diagonaux des classes de $f_n$ et de $\gamma_n$
modulo $\gothq$ donnent des solutions pour l'\'{e}quation de
Maurer-Cartan pour chaque composante $F_a\subset E$
individuellement. Mais $\Ext^2(F_a,F_a)_0=0$ et ainsi il n'y a pas
d'obstruction pour une telle d\'{e}formation. On conclut que $f\equiv 0
\mod \gothq$.
\end{proof}

%%%%%%%%%%%%%%%%%%%%%%%%%%%%%%%%%%%%%%%%%%%%%%%%%%%%%%%%%%%%%%%%%%%%%%%%%%
%%%
%%%   Referenzliste
%%%
%%%%%%%%%%%%%%%%%%%%%%%%%%%%%%%%%%%%%%%%%%%%%%%%%%%%%%%%%%%%%%%%%%%%%%%%%%
\bibliographystyle{plain}

\parindent0mm

\end{document}